\documentclass[12pt]{article}
\usepackage{amsfonts,amsmath}

                    \def\cF{{\cal F}}

                    \def\cR{{\cal R}}
                    \def\cU{{\cal U}}

\newcommand{\CC}{{\mathbb C}}

\newcommand{\ZZ}{{\mathbb Z}}

\newcommand{\nn}{\nonumber\\ }

\newcommand{\Rtriang}[9]{
  {\left(
    \begin{array}{cccc}
      1 & \ #1\  & \ #2\  & \ #3 \ \cr
      & #4 & #5 & \ #6\  \cr
      & & \ #7\  & \ #8\  \cr
      & & & #9
    \end{array}
  \right)
  }}

\normalbaselines 
\begin{document}
\pagestyle{empty}

\begin{center}
  \textsf{\LARGE 
  On Combined Standard-Nonstandard\\[3mm]
 or Hybrid $(q,h)$-Deformations
  }

\vspace{7mm}

{\large B.L.~Aneva$^{a,}$\footnote{blan@inrne.bas.bg}, 
 ~D.~Arnaudon$^{b,}$\footnote{Daniel.Arnaudon@lapp.in2p3.fr}, 
 ~A.~Chakrabarti$^{c,}$\footnote{chakra@cpht.polytechnique.fr},\\[2mm] 
 V.K.~Dobrev$^{d,a,}$\footnote{vladimir.dobrev@unn.ac.uk,dobrev@inrne.bas.bg} 
~and~ S.G.~Mihov$^{a,}$\footnote{smikhov@inrne.bas.bg}}

\vspace{5mm}

\emph{$^a$ Institute of Nuclear Research and Nuclear Energy} 
\\
\emph{Bulgarian Academy of Sciences}
\\
\emph{72 Tsarigradsko Chaussee, 1784 Sofia, Bulgaria}
\footnote{permanent address for V.K.D.}
\\
\vspace{3mm}
\emph{$^b$ Laboratoire d'Annecy-le-Vieux de Physique Th{\'e}orique LAPTH}
\\
\emph{CNRS, UMR 5108, associ{\'e}e {\`a} l'Universit{\'e} de Savoie}
\\
\emph{LAPP, BP 110, F-74941 Annecy-le-Vieux Cedex, France}
\\
\vspace{3mm}
\emph{$^c$ Centre de Physique Th{\'e}orique, CNRS UMR 7644}
\\
\emph{Ecole Polytechnique, 91128 Palaiseau Cedex, France.}
\\
\vspace{3mm}
\emph{$^d$ School of Computing and Mathematics}
\\
\emph{University of Northumbria}
\\
\emph{Ellison Place, Newcastle upon Tyne, NE1 8ST, UK}

\end{center}

\vspace{.8 cm}

\begin{abstract}
  Combined $(q,h)$-deformations proposed by Kupershmidt and
  Balles\-teros--Herranz--Parashar are studied. In each case a
  transformation is shown to lead to an equivalent, standard 
  $q$-deformation. We briefly indicate that appropriate singular
  limits of the same type of transformations can however lead from
  standard biparametric $(p,q)$-deformations to non-hybrid but
  biparametric nonstandard $(g,h)$ ones. Finally a case of
  hybrid $(q,h)$-deformation is recalled, related to the superalgebra 
  $GL(1|1)$. 
\end{abstract}

\vfill

\rightline{INRNE-TH-00-02}
\rightline{LAPTH-800/00}
\rightline{S 080.0600}
\rightline{UNN-SCM-M-00-04}
\rightline{math.QA/0006206}
\rightline{June 2000}

\newpage
\pagestyle{plain}
\setcounter{page}{1}

\section{Introduction}
\label{sect:intro}
\setcounter{equation}{0}
Not only the standard $q$-deformation but also the nonstandard
(Jordanian) $h$-deform\-ation of $GL(2)$ can be considered
to be well-known. In each of these domains biparametric
generalizations, $(p,q)$ and $(g,h)$ respectively, have been studied
by number of authors. A large number of previous sources are cited in
\cite{Kup} and \cite{BHP}. The dual quantum algebras of $GL_{pq}$ 
and $GL_{gh}$ were found in \cite{DobrevJMP33} and
\cite{ADMjpa30}, respectively. 
Here we are concerned with certain proposals for combining these two
distinct types into $(q,h)$-deformations. They will often be denoted
as hybrid ones. In particular, we analyze the results of Kupershmidt
\cite{Kup} and of Ballesteros--Herranz--Parashar \cite{BHP}. 
In each case, we show that a well-defined transformation eliminates
$h$ leaving a standard  $q$-deformation. 
%% (One may use the
%% terminology ``dehybridization'' in this context.) 
This transformation
is not a arbitrary twist, but a straightforward
similarity relation performed by a tensor square of an operator. 
This will be demonstrated explicitly in Sect. \ref{sect:Kup} and
\ref{sect:BHP}. 

In  Sect. \ref{sect:Kup} we start, in fact, with the $3$-parameter
$(q,h,h')$ deformation of \cite{Kup}. Already at this level we are
able to construct a similarity transformation reducing the formalism
to a 1-parameter deformation. The surviving single parameter $q'$ is
expressed explicitly as a function of $(q,h,h')$. For the case
particularly advocated in \cite{Kup}, namely $h'=0$, one has simply
$q'=q$. 

In Sect. \ref{sect:BHP} we start by transforming the $R_{q,h}$ matrix
of \cite{BHP} to  $R_q$. Then we show how their relevant results
can be much better understood in the context of an explicitly
presented, ``coalgebra conserving'' map. This illuminates several
aspects and goes beyond the case of $4\times 4$ matrices. 

In Sect. \ref{sect:singular} we add some comments on maps 
and singular limits of transformations. Their nontrivial consequences
\cite{ACCS,ACCtowards,CCgnf,CQjord} are indicated by
citing appropriate references. The passage from a standard
biparametric $(p,q)$-deformation to a nonstandard $(g,h)$ one is
presented in this context. 

In Sect. \ref{sect:authentic} we arrive (at last) to a 
hybrid $(q,h)$ deformation where $h$ cannot be transformed away. 
This turns out to be a hybrid 
deformation of the superalgebra $GL(1|1)$,
already studied in \cite{FHRclass}. It 
is located as the case $R_{H1.2}$ in the classification of $4\times 4$
$R$-matrices in \cite{Hietarinta}, which we briefly recall. 

Finally we would like to come back to Sections \ref{sect:Kup} and
\ref{sect:BHP}. Instead of briefly stating the equivalence
$(q,h)\longrightarrow (q)$, we have chosen to present our elementary
analysis explicitly and in some detail. We consider this worthwhile
for dissipating some confusions. Several authors have presented
attractive looking hybrid deformations without noticing disguised
equivalences. We ourselves devoted time and effort to their study
before reducing them to usual deformations. We hope that our analysis
will create a more acute awareness of traps in this domain.

\section{Kupershmidt's $(q,h,h')$ and $(q,h)$ deformations}
\label{sect:Kup}
\setcounter{equation}{0}
We start by noting that the group relations given by the set of
equations (5) of \cite{Kup} can be written as 
\begin{eqnarray}
  \label{eq:Kup5}
  ca &=& ac \;, \qquad \qquad \qquad bd\ =\ db \nonumber\\
  cb &=& qbc -hac -h'db \nonumber\\
  ad &=& da +(q-1)bc -hac -h'db \nonumber\\
  qba &=& ab +ha^2 + h'b^2 - h(da-bc) \nonumber\\
  cd &=& qdc -hc^2 -h'd^2 +h'(da-bc) \;.
\end{eqnarray}
{}From the second and the third equations of (\ref{eq:Kup5}) one
obtains 
\begin{equation}
  \label{eq:Kup+}
  ad - qbc + hac + h'bd = da-bc = ad-cb \;.
\end{equation}
Substituting from (\ref{eq:Kup+}) the l.h.s. for $(da-bc)$ in the
fourth and fifth equations of (\ref{eq:Kup5})
one gets back exactly (5) of \cite{Kup}. Compared to his original
version ours has the following advantages
\begin{itemize}
\item Adopting the ordering 
  \begin{equation}
    \label{eq:dabc}
    d>a>b>c   
  \end{equation}
  all the terms in increasing
  order ($ca$, $bd$,...) are on the l.h.s. of (\ref{eq:Kup5}), whereas
  the square terms and terms in decreasing order are one the
  r.h.s. This solves the ordering problem encountered in \cite{Kup}
  when $h'$ was taken different from 0.
\item The roles of the parameters $h,h'$ are now more simple and
  symmetrical. The terms bilinear in them (like $hh'bd$ and ${h'}^2bd$)
  do not appear in (\ref{eq:Kup5}). The corresponding complementary
  (upper and lower triangular) linear contributions of $h$ and $h'$ in
  the $R$-matrix to be presented below correspond directly to this
  feature (the possibility of linearizing the contributions in
  (\ref{eq:Kup5}). 
\item The simpler form of (\ref{eq:Kup5}) facilitates the construction
  of the $R$-matrix form of the RTT relations. This, in turn,
  facilitates the construction of the explicit similarity 
  transformation leading to a $1$-parametric equivalent deformation
  \begin{equation}
    \label{eq:qhh'2q'}
    (q,h,h') \longleftrightarrow q'
  \end{equation}
  where $q'$ is a specific function of $(q,h,h')$ to be presented
  below. We could have derived the same final results using the more
  complicated version of (\ref{eq:Kup5}) in \cite{Kup}. But
  (\ref{eq:Kup5}) is preferable. 
\end{itemize}

\subsubsection*{Solution of the RTT constraints}

Let 
\begin{equation}
  \label{eq:T}
  T = \left(
    \begin{array}{cc}
      a & b\\
      c & d
    \end{array}
  \right)\;,
  \qquad\qquad 
  T_1 = T\otimes 1, \qquad T_2 = 1\otimes T \;,
\end{equation}
where $a,b,c,d$ satisfy (\ref{eq:Kup5}) and let $R$ satisfy 
\begin{equation}
  \label{eq:RTT}
  RT_1 T_2 = T_2 T_1 R \;.
\end{equation}
Using (\ref{eq:Kup5}) systematically, one obtains a solution involving
an arbitrary parameter $\kappa$. It is 
\begin{equation}
  \label{eq:Rhh'}
  R = \left(
    \begin{array}{cccc}
      1 & -h\kappa & h\kappa & 0 \\
      0 & q\kappa & (1-q\kappa) & 0 \\
      0 & (1-\kappa) & \kappa & 0 \\
      0 & -h'\kappa & h'\kappa & 1   
    \end{array}
    \right) \;.
\end{equation}
This does \emph{not} satisfy the Yang--Baxter (YB) relations for all
values of $\kappa$. In fact, in order that $R$ satisfies YB the parameter 
$\kappa$ must satisfy the following quadratic equation:
\begin{equation}
  \label{eq:KK}
\kappa^2 (q+h h') - \kappa (q+1) +1 =0
\end{equation}
(Thus, for example, the particularly simple form for $\kappa=1$
does not satisfy YB unless $h h'=0$.) 
The presence of $\kappa$ at this stage permits the
existence of two solutions satisfying YB constraints and related
through $R\longleftrightarrow R^{-1}_{21}$. 

We find it convenient for the 
similarity transformation to be introduced below  
to write down the two solutions of (\ref{eq:KK}) 
in the following manner: 
\begin{equation}
  \label{eq:K}
  \kappa=\kappa_1 = (1+\eta^{-1}h)^{-1} \ ,
  \qquad
  \kappa=\kappa_2 = (1+\eta h')^{-1} 
\end{equation}
where the parameter $\eta$ satisfies the quadratic 
\begin{equation}
  \label{eq:eta}
  \eta^{-1} h + \eta h' = q-1
\end{equation}

Our equation (\ref{eq:eta}) is the same as the one used in \cite{Kup}
eq. (15) to eliminate $h'$ at the level of the vector basis of the
Poisson bracket algebra. But the r{\^o}le of our $\eta$ (corresponding
to $t$ in \cite{Kup}) is different. We continue to allow $h'$ to be
arbitrary and finally use a similarity transformation to arrive at an
equivalent 1-parameter deformation. 

\subsubsection*{Similarity transformation to a 1-parameter
  deformation}  

Define 
\begin{equation}
  \label{eq:G}
  G = \left(
    \begin{array}{cc}
      1 &\eta \\
      \zeta & (1+\eta\zeta)
    \end{array}
  \right)
  \qquad\qquad
  G^{-1} = \left(
    \begin{array}{cc}
      (1+\eta\zeta) & -\eta \\
      -\zeta & 1
    \end{array}
  \right)
\end{equation}
where the parameters $\eta,\zeta$ are given by 
\begin{equation}
  \label{eq:hh'}
  h=(q-1)\frac{\eta(1+\eta\zeta)}{1+2\eta\zeta} \;,
  \qquad\qquad
  h'=(q-1)\frac{\eta}{1+2\eta\zeta} \;,
\end{equation}
or 
\begin{equation}
  \label{eq:etazeta}
\begin{array}{l}
  \eta = (2h')^{-1} ((q-1)\pm \lambda) \,,
  \qquad
  \eta^{-1} = (2h)^{-1} ((q-1)\mp \lambda) \\ \\ 
   \zeta = \mp h' \lambda^{-1} \,, 
  \qquad 
  \lambda=\sqrt{(q-1)^2-4hh'} 
\end{array}
\end{equation}
(For $\lambda=0$, $\zeta$ diverges. This point should be approached as a
limit after transforming.)

On obtains, after subtle simplifications
\begin{equation}
  \label{eq:R'}
  R'=(G\otimes G)R(G^{-1}\otimes G^{-1})
  = \left(
    \begin{array}{cccc}
      1 & 0 & 0 & 0 \\
      0 & \kappa/\kappa_1 & (1-\kappa/\kappa_1) & 0 \\
      0 & (1-\kappa/\kappa_2) & \kappa/\kappa_2 & 0 \\
      0 & 0 & 0 & 1   
    \end{array}
    \right) \;.
\end{equation}
where $\kappa_1,\kappa_2$ are given by (\ref{eq:K}). (In the light of the
remarks following (\ref{eq:etazeta}), note that $\zeta$ does not appear in
$R'$, only $\eta$ through $\kappa_1$ and $\kappa_2$. Reality restrictions are
discussed at the end of this section.) 

Now the statement leading to (\ref{eq:K}) is evident: for $\kappa=\kappa_1$ one
gets 
\begin{equation}
  \label{eq:R'1}
  R'
  = \left(
    \begin{array}{cccc}
      1 & 0 & 0 & 0 \\
      0 & 1 & 0 & 0 \\
      0 & (1-\kappa_1/\kappa_2) & \kappa_1/\kappa_2 & 0 \\
      0 & 0 & 0 & 1   
    \end{array}
    \right) \;,
\end{equation}
the standard lower triangular form of the YB solution for the single
parameter $\kappa_1/\kappa_2$. Similarly, for $\kappa=\kappa_2$, one gets 
\begin{equation}
  \label{eq:R'2}
  R'
  = \left(
    \begin{array}{cccc}
      1 & 0 & 0 & 0 \\
      0 & \kappa_2/\kappa_1 & (1-\kappa_2/\kappa_1) & 0 \\
      0 & 0 & 1 & 0 \\
      0 & 0 & 0 & 1   
    \end{array}
    \right) \;,
\end{equation}
the standard upper-triangular YB solution for the single parameter
$\kappa_2/\kappa_1$. In fact (\ref{eq:R'1}) and (\ref{eq:R'2}) are related as
the pair $R'$, $(R'_{21})^{-1}$.

Let us now define 
\begin{equation}
  \label{eq:T'}
    T' = \left(
    \begin{array}{cc}
      a' & b'\\
      c' & d'
    \end{array}
  \right)\;,
\end{equation}
such that it satisfies
\begin{equation}
  \label{eq:R'T'T'}
  R'T'_1 T'_2 = T'_2 T'_1 R' \;.
\end{equation}
The group relations turn out to be independent of the parameter $\kappa$ of
(\ref{eq:R'}) (just as (\ref{eq:Kup5}) are independent of $\kappa$ in
(\ref{eq:Rhh'})). One obtains
\begin{eqnarray}
  \label{eq:RTTsl2q1}
  c'a' &=& a'c', \qquad 
  b'd' \ =\  d'b', \qquad
  c'b' \ =\  q'b'c', \nonumber\\
  a'd' &=& d'a'+(q'-1) b'c',  \qquad
  q'b'a' \ =\  a'b', \qquad
  c'd' \ =\  q'd'c',
\end{eqnarray}
where we have set 
\begin{equation}
  \label{eq:q'}
  q' = \frac{\kappa_2}{\kappa_1} = \frac{1+\eta^{-1}h}{1+\eta h'}
\end{equation}
with $\eta^{-1}h+\eta h'=q-1$.
Note that 
\begin{eqnarray}
  &&\mbox{for $h'=0$,}\qquad q'=q \nonumber\\
  &&\mbox{for $h=0$,}\qquad q'=q^{-1} 
\end{eqnarray}
A more complete discussion of the domains of $q'$ follows below.

Setting $h=h'=0$ in (\ref{eq:Kup5}) and adding primes one obtains
(\ref{eq:RTTsl2q1}). 

The relations (\ref{eq:RTTsl2q1}) can also be obtained from
(\ref{eq:Kup5}) by transforming with $G$. We have preferred to
construct the corresponding $R$-matrices first from their intrinsic
interest and also for elucidating the significance of the free
parameter $\kappa$ arising at the RTT level before imposing the YB
constraints. 

In \cite{Kup} after setting $h'=0$ the $(q,h)$ deformed system is
reformulated using a certain ordering. Our preceding study contains
the $R$-matrix for this case as the particular one obtained by setting
in (\ref{eq:Rhh'})
\begin{eqnarray}
  \kappa &=& \kappa_2 = (1+\eta h')^{-1} = 1 \qquad \mbox{(for $h'=0$)} \\
  \kappa_1 &=& (1+\eta^{-1}h)^{-1} = q^{-1} \nonumber
\end{eqnarray}
giving 
\begin{equation}
  \label{eq:Rhybrid}
  R'
  = \left(
    \begin{array}{cccc}
      1 & -h & h & 0 \\
      0 & q & (1-q) & 0 \\
      0 & 0 & 1 & 0 \\
      0 & 0 & 0 & 1   
    \end{array}
    \right) \;.
\end{equation}
The transformation to 1-parameter form is now by 
\begin{equation}
  \label{eq:Ghybrid}
    G = \left(
    \begin{array}{cc}
      1 &\eta \\
      0 & 1
    \end{array}
  \right)
  \qquad \mbox{with $\zeta=0$, $\eta=\frac{h}{q-1}$}
\end{equation}
in (\ref{eq:G}). The upper triangular form (\ref{eq:R'2}) now becomes 
\begin{equation}
  \label{eq:Rsl2q1bis}
  R' 
  = \left(
    \begin{array}{cccc}
      1 & 0 & 0 & 0 \\
      0 & q & (1-q) & 0 \\
      0 & 0 & 1 & 0 \\
      0 & 0 & 0 & 1   
    \end{array}
    \right) \;.
\end{equation}
The transformed $(a,b,c,d)$, namely $(a',b',c',d')$ now satisfies
(\ref{eq:RTTsl2q1}) with $q'=q$. 

The results for $h'=0$ can be obtained directly by starting with a $(q,h)$
system given by (16) of \cite{Kup} with one necessary correction. The
group relations should be written as 
\begin{eqnarray}
  \label{eq:Kup16}
  ba &=& q^{-1} a(b+ha) - q^{-1} hda + q^{-2} hc (b+ha) \nonumber \;, \\
  bd &=& db \nonumber \;, \\
  bc &=& q^{-1} c(b+ha) \nonumber \;, \\
  ad &=& da+ (1-q^{-1})cb - q^{-1} hca \nonumber \;, \\
  ac &=& ca \nonumber \;, \\
  dc &=& q^{-1} c(d+hc) \;.
\end{eqnarray}
(The last two terms of the first equation has each an extra factor
$q^{-1}$ as compared to (16) of \cite{Kup}. Setting $h'=0$ in (5) of
\cite{Kup} and reordering one indeed gets our version. Starting
directly with (\ref{eq:Kup16}), the RTT relations and the YB
constraints can indeed be shown to lead to (\ref{eq:Rhybrid}). Then
(\ref{eq:Ghybrid}) eliminates $h$ leading to (\ref{eq:Rsl2q1bis}) and to
(\ref{eq:RTTsl2q1}) with $q'=q$. The case $h=0$ can be treated quite
analogously setting (taking the lower sign in (\ref{eq:etazeta}) for
$\eta,\zeta$) 
\begin{equation}
  \eta=0, \qquad 
    G = \left(
    \begin{array}{cc}
      1 & 0 \\
      \zeta & 1
    \end{array}
  \right)
  , \qquad
  \zeta = \frac{h'}{q-1} \;.
\end{equation}

In \cite{Kup} $h'$ was eliminated at the stage of Poisson brackets and
$h$ was retained. It was assumed that one thus obtains an authentic
2-parameter $(q,h)$-deformation. We have shown that this is not the
case. Our transformation for the original $(q,h,h')$ case shows that,
from the start, one has always been dealing with a heavily disguised
$q$-deformation. This statement should however be qualified by taking
a closer look at different domains of the parameter space of
$(q,h,h')$. We consider below real values of $(q,h,h')$. 

For both $(h,h')\neq 0$ (the cases $h'=0$ and $h=0$, with $q'=q$ and
$q'=q^{-1}$ respectively, 
can be considered  simply and analogously), from (\ref{eq:q'}),
\begin{equation}
  \label{eq:q'bis}
  q' = \frac{q+1\pm \sqrt{(q-1)^2-4hh'}}{q+1\mp \sqrt{(q-1)^2-4hh'}} 
\end{equation}

Apart from the very special case 
\begin{equation}
  \label{eq:q=-1}
  q=-1 \qquad\qquad q'=-1
\end{equation}
we note that more generally
\emph{i)} for any $q$ and $hh'<0$
and \emph{ii)} for $(q-1)^2>4hh'>0$, $q'$ is always real. We consider
this as the generic case. 

Another very special case is (for $hh'>0$) $|q-1|=2\sqrt{hh'}$ when
$q'=1$. See the remark following (\ref{eq:etazeta}) concerning this
singular point. Here one has a classical solution with commuting
$a',b',c',d'$~!

For $hh'>0$ and   $|q-1|<2\sqrt{hh'}$, $q'=e^{\pm i\delta}$, a complex
phase. Here $q'$ can even be a root of unity. Thus starting from a
complex deformation, one can obtain by the transformation with $G$ an
equivalent deformation with 3 real parameters related through
(\ref{eq:q'bis}).

\section{The Ballesteros--Herranz--Parashar (BHP) case}
\label{sect:BHP}
\setcounter{equation}{0}

\subsection{Transformation to $R_q$}
The BHP two-parametric deformation (Sect. 4 of \cite{BHP}) leads to
the $R$-matrix
\begin{equation}
  \label{eq:RBHP}
  R_{q,h} 
  = \left(
    \begin{array}{cccc}
      1 & h & -qh & h^2 \\
      0 & q & 1-q^2 & qh \\
      0 & 0 & q & -h \\
      0 & 0 & 0 & 1   
    \end{array}
    \right) \;.
\end{equation}
The authors present it as a superposition of standard ($q$) and
nonstandard ($h$) deformations. It has indeed the attractive property
that for $h=0$ and $q=1$ respectively one obtains the standard
$R_q$ and the nonstandard $R_h$ matrices. 

Now consider a similarity transformation of $R_{q,h}$ by $M\otimes M$
where 
\begin{equation}
  \label{eq:M}
    M = \left(
    \begin{array}{cc}
      x & y \\
      0 & 1
    \end{array}
  \right)
  , \qquad
    M^{-1} = x^{-1} \left(
    \begin{array}{cc}
      1 & -y \\
      0 & x
    \end{array}
  \right)
\end{equation}
with $y=\frac{\displaystyle h}{\displaystyle q-1}$ ($h\neq 0$, $q\neq
1$) and $x$ is an arbitrary nonzero parameter. In
the notation of \cite{BHP}
\begin{equation}
  \label{eq:y}
  y=\frac{a_+}{2a} \qquad (q=e^a)
\end{equation}
One obtains, independently of the choice of $x$,
\begin{equation}
  \label{eq:Rq}
  (M^{-1} \otimes M^{-1}) R_{q,h} (M \otimes M)
  =
   \left(
    \begin{array}{cccc}
      1 & 0 & 0 & 0 \\
      0 & q & 1-q^2 & 0 \\
      0 & 0 & q & 0 \\
      0 & 0 & 0 & 1   
    \end{array}
    \right) 
    = R_q \;.
\end{equation}
Thus it is seen that $h$ can be transformed away. 

\subsection{A coalgebra preserving map}
The significance of this equivalence (and that of the related Hopf
algebraic results of \cite{BHP}) are better understood in the context
of a simple class of coalgebra preserving maps, presented below. (The
content of the mapping can be considered, in a certain sense, to be
trivial. Elucidating this aspect is precisely our purpose.) These maps
can be generalized to higher dimensional algebras. But we here
consider only $\cU_q(gl(2))$.

One starts with the standard $\cU_q(sl(2))$ algebra
\begin{eqnarray}
  \label{eq:Uqsl2}
  {}[J_0,J_\pm] &=& \pm 2J_\pm \qquad\qquad 
  (q^{J_0}J_\pm = q^{\pm 2}J_\pm q^{J_0}) \nonumber\\
  {}[J_+,J_-] &=& \frac{q^{J_0}-q^{-J_0}}{q-q^{-1}} \equiv [J_0]
\end{eqnarray}
with the coalgebra structure
\begin{eqnarray}
  \label{eq:Uqsl2co}
  \Delta (J_0) &=& J_0\otimes 1+ 1 \otimes J_0 \qquad\qquad
  (\Delta(q^{J_0}) = q^{J_0} \otimes q^{J_0})
  \nonumber\\
  \Delta (J_\pm) &=& q^{J_0/2} \otimes J_{\pm} + 
  J_{\pm} \otimes q^{-J_0/2}
  \nonumber\\
  S(q^{J_0/2}) &=& q^{-J_0/2} \nonumber\\
  S(J_{\pm}) &=& - q^{-J_0/2} J_{\pm} q^{J_0/2} \nonumber\\
  \epsilon (q^{J_0/2}) &=& 1, \qquad\qquad \epsilon(J_{\pm})=0
\end{eqnarray}
(Note that for the antipode $S$, we have not replaced 
$q^{-J_0/2}J_\pm q^{J_0/2}$ by $q^{\mp 1}J_\pm$. This last form uses
(\ref{eq:Uqsl2}) which will be modified by the map, keeping the
structure (\ref{eq:Uqsl2co}) intact in terms of the new generators.)
 
Next one sets
\begin{eqnarray}
  \label{eq:map}
  J'_0 &=& J_0 \nonumber\\
  J'_+ &=& J_+ 
  %%a_1 J_+ + a_2(q^{J_0/2} - q^{-J_0/2}) + a_3 J_- 
  \nonumber\\
  J'_- &=& b_1 J_- + b_2(q^{J_0/2} - q^{-J_0/2}) + b_3 J_+ \;.
\end{eqnarray}
Choosing (with $q=e^a$) 
\begin{eqnarray}
  \label{eq:a1b1}
%%  a_1&=&1, \qquad a_2\ =\ 0, \qquad a_3\ =\ 0, \nonumber\\
  b_1 &=& \frac{\sinh(a)}{a}, \qquad b_2\ =\ - \frac{a_+}{2a^2}
  \qquad b_3\ =\ -\frac{a_+^2}{4a^2}
\end{eqnarray}
one obtains the case (4.6) of \cite{BHP} (our $(J'_0,J'_+,J'_-)$
corresponding to their $(J'_3,J_+,J_-)$).
Indeed
\begin{eqnarray}
  \label{eq:relcomBHP}
  {}[J'_0,J'_+] &=&  2J'_+ \nonumber\\
  {}[J'_0,J'_-] &=&  - 2J'_- - \frac{a_+}{a} \frac{\sinh(aJ'_0/2)}{a/2} 
  - \frac{a_+^2}{a^2} J'_+ \nonumber\\
  {}[J'_+,J'_-] &=& \frac{\sinh(aJ'_0)}{a} + \frac{a_+}{a}
  \frac{e^a-1}{2a} (e^{-aJ'_0/2} J'_+ + J'_+ e^{aJ'_0/2}) 
\end{eqnarray}
Moreover, the coalgebra structure induced by this map has the same
expression in terms of the new generators, i.e.
\begin{eqnarray}
  \label{eq:DeltaJ'}
  \Delta (J'_0) &=& J'_0\otimes 1+ 1 \otimes J'_0 \nonumber\\
    \Delta (J'_\pm) &=& q^{J'_0/2} \otimes J'_{\pm} + 
    J'_{\pm} \otimes q^{-J'_0/2}
\end{eqnarray}
and so on. This is achieved under the single condition that the
coefficients of $q^{\pm J_0}$ in (\ref{eq:map}) are opposite. 
(Note that the latter statement would be true also if, in addition, we
write 
  $J'_+ = a_1 J_+ + a_2(q^{J_0/2} - q^{-J_0/2}) + a_3 J_- $.)
The whole Hopf algebra described by $(J'_0,J'_+,J'_-)$ then
reproduces that of BHP. 
This means that the BHP Hopf algebra is equivalent to
$\cU_q(gl(2))$.

\section{Comments on singular limits of transformations}
\label{sect:singular}
\setcounter{equation}{0}

Up to now we have been considering regular, invertible transformations
making evident the trivial nature of the passage 
\begin{equation}
  \label{eq:passage}
  R_q \longleftrightarrow R_{q,h} \;.
\end{equation}
%% (Note that we do not have 
%% $  R_h \longleftrightarrow R_{q,h} $ consistently 
%% with the comments at the end of section \ref{sect:BHP}.) 
The map (\ref{eq:map}) is consistent with this
due to the conservation of the structure of the coalgebra. 
When a map has to be followed by a twist \cite{DrinQuasi,Resh:lmp20}
to arrive at a sought for
coalgebra, the situation can be of interest. It has been shown
elsewhere \cite{ACCS} how the universal $\cR_h$ matrix 
(introduced first in \cite{Ogiev})
can be obtained, through a twist, starting from the trivial classical one
($\cU_h(sl(2))$ is a triangular Hopf algebra). Another interesting
possibility is the use of a transformation singular in the limit
$q\rightarrow 1$ but in such a specific fashion that ($G(q,h)$ being
singular at $q=1$),
\begin{equation}
  \label{eq:singularmap}
  \left(
    G(q,h)^{-1}\otimes G(q,h)^{-1} R_q G(q,h)\otimes G(q,h)
  \right) \Big|_{q=1}
  = R_h \;.
\end{equation}
In contrast to (\ref{eq:passage}) this passage is non-invertible and
the end-product is not a hybrid $R_{q,h}$ but a nonstandard
$R_h$. 
This can be considered as an operator equation between universal
R-matrices and $G$ given by (as shown in \cite{ACCtowards}): 
%%It has been shown elsewhere \cite{ACCtowards} that a correct
%%prescription is ($J_+$ satisfying (\ref{eq:Uqsl2}))
\begin{eqnarray}
  \label{eq:Gqh}
  G(q,h) &=& E_q(\eta J_+) \qquad \qquad \mbox{with}\qquad
  \eta=\frac{h}{q-1} \\
  \mbox{and} \qquad\qquad 
  E_q(x) &=& \sum_{n=0}^\infty \frac{x^n}{[n]!}, \qquad
  [n]\equiv \frac{q^n-q^{-n}}{q-q^{-1}}
\end{eqnarray}
Alternatively (4.2) can be regarded as a matrix equation implementing
$j_1 \otimes j_2$ representations. 
This technique can be generalized to $GL(N)_q$ \cite{ACCtowards} and
also to obtain nonstandard quasi-Hopf algebras \cite{CCgnf}.

Note that in (\ref{eq:singularmap}) $\eta$ has the same form as in
(\ref{eq:G}) or (\ref{eq:M}). But the crucial difference is that one
takes the limit $q\longrightarrow 1$.

In the above-mentioned references universal $\cR$ matrices have been
studied but only for $R_q$ and $R_h$. Here, in conclusion, we indicate
how one can treat the biparametric case involving $R_{p,q}$ and
$R_{g,h}$. (Note that no hybrid deformation is involved here.)
We restrict the study to the fundamental case of $4\times 4$
matrices. 

We start with 
\begin{equation}
  \label{eq:Rpq}
  R_{pq} = 
  \left(
    \begin{array}{cccc}
      p & 0 & 0 & 0 \\
      0 & pq & p-q & 0 \\
      0 & 0 & 1 & 0 \\
      0 & 0 & 0 & p   
    \end{array}
    \right) 
\end{equation}
and define (with $j_1\otimes j_2 = \frac12 \otimes \frac12$ in
(\ref{eq:Gqh})) 
\begin{equation}
    G = \left(
    \begin{array}{cc}
      1 & \eta \\
      0 & 1
    \end{array}
  \right)
\end{equation}
but now $\eta$ will be chosen differently. One obtains
\begin{equation}
  \label{eq:Rpqeta}
  (G^{-1}\otimes G^{-1}) R_{p,q} (G\otimes G) = 
  \left(
    \begin{array}{cccc}
      p & p(1-q)\eta & (q-1)\eta & (1-p)(q-1)\eta^2 \\
      0 & pq & p-q & q(p-1)\eta \\
      0 & 0 & 1 & (1-p)\eta \\
      0 & 0 & 0 & p   
    \end{array}
    \right) 
    \;.
\end{equation}
Now let $q\rightarrow 1$ and $p\rightarrow 1$ in such a fashion that 
\begin{equation}
  \label{eq:lambda}
  \left( \frac{q-1}{p-1} \right)^{1/2} = \lambda = \mbox{constant.}
\end{equation}
Set
\begin{equation}
  \label{eq:etalambda}
  \eta = \frac{\eta_0}{((p-1)(q-1))^{1/2}}, \qquad\qquad
  \lambda\eta_0=h, \qquad\qquad \lambda^{-1}\eta_0 = -g \;.
\end{equation}
Then 
\begin{equation}
  \label{eq:Rgh}
  \left.\left(
    (G^{-1}\otimes G^{-1}) R_{p,q} (G\otimes G) 
  \right)\right|_{(p\rightarrow 1,q\rightarrow 1)}
  = \left(
    \begin{array}{cccc}
      1 & -h & h & gh \\
      0 & 1 & 0 & -g \\
      0 & 0 & 1 & g \\
      0 & 0 & 0 & 1   
    \end{array}
    \right) 
    \;.
\end{equation}
Thus we obtained the 2-parametric nonstandard $R$-matrix. For $g=h$
and $g=0$, we obtain the two known forms of $R_h$. The 2-parametric
universal coloured, nonstandard $R$-matrix for deformed $gl(2)$ is
obtained in Sect. 3 of \cite{CQjord} implementing a twist. Here we
presented the $4\times 4$ case to show how $\eta$ is quite simply
modified from (\ref{eq:Gqh}) to (\ref{eq:etalambda}) as one passes
from the 1- to the 2-parametric case.

\section{An authentic hybrid $(q,h)$ deformation: $GL_{q,h}(1|1)$}
\label{sect:authentic}
\setcounter{equation}{0}
In Sect. \ref{sect:Kup} and Sect. \ref{sect:BHP} we have shown that the 
hybrid $(q,h)$ deformations in \cite{Kup} and \cite{BHP} are in fact
disguised $q$-ones.  In the search of hybrid deformation we also check 
with the classification of $4\times 4$ $R$-matrices in \cite{Hietarinta}. 
There we find seven triangular cases: 

\begin{equation}
  R_{S2,1} = \Rtriang{0}{0}{0}
  {p}{1-pq}{0}
  {q}{0}
  {1}
  \end{equation}

\begin{equation}
  R_{S2,2} = \Rtriang{0}{0}{0}
  {p}{1-pq}{0}
  {q}{0}
  {-pq}
\end{equation}

\begin{equation}
  (R_{H1,3})_{\vert_{k=1,p=-h,q=-g}}\ =\ \Rtriang{-h}{h}{gh}
  {1}{0}{-g}
  {1}{g}
  {1}
\end{equation}

\begin{equation}
  (R_{H2,3})_{\vert_{k=1},p=x_1,q=x_2,s=x_3}\ =\ \Rtriang{x_1}{x_2}{x_3}
  {1}{0}{x_2}
  {1}{x_1}
  {1}
\end{equation}

\begin{equation}
  R_{S0,1} = \Rtriang{0}{0}{1}
  {1}{0}{0}
  {1}{0}
  {1}
\end{equation}

\begin{equation}
  R_{S0,2} = \Rtriang{0}{0}{1}
  {-1}{0}{0}
  {-1}{0}
  {1}
\end{equation}

\begin{equation}
\label{eq:Rqh}
 R_{q,h}\ =\  (R_{H1,2})_{\vert_{p=1,k=h}}\ =\ \Rtriang{0}{0}{h}
  {1}{1-q}{0}
  {q}{0}
  {-q}
\end{equation}

Note that in \cite{Hietarinta} the $R$-matrices are given in 
two versions:  homogeneous $R_{H...}$ and scaled $R_{S...}$. 
The scaled versions are simpler but in some cases, in order not 
to lose some symmetry among the parameters 
we use the homogeneous versions with only 
an overall rescaling. 

The case $R_{S2,1}$ is the 2-parameter $p,q$ deformation, 
$GL_{pq}(2)$ the dual of which is given in \cite{DobrevJMP33}. 
The case $R_{S2,2}$ is a superalgebra - 
the known $p,q$ deformation of $GL_{pq}(1|1)$, 
the dual of which is given in \cite{HiRi,DaWa,BuTo}.

\textit{Remark:} Note that here (and below for $R_{q,h}$)
we consider ordinary, not graded, $R$-matrices. 
The results can be translated to the graded formalism. 
There is a one-to-one correspondence between the results 
obtained through the two approaches (which was noticed first 
for solutions of YBE and graded YBE in \cite{KulSkly}). 
This correspondence may be given also 
through ``transmutation'' in the 
sense of \cite{Majidbook}. This aspect is considered in the context of
$sl(1|2)$ in \cite{ACFsl12qs}. 
In our paper the superalgebraic aspect becomes
evident after implementation of the RTT formalism.

In the third case we have written the homogeneous version 
$R_{H1,3}$ of \cite{Hietarinta} 
with $k=1$ and renamed parameters. This seems the natural scaling 
(and not the $R_{S1,3}$). The result is indeed with 
2 parameters, i.e., this is the 2-parameter Jordanian 
$GL_{gh}(2)$. The dual was found in \cite{ADMjpa30}.

In the fourth case we have written the homogeneous version 
$R_{H2,3}$ of \cite{Hietarinta} with $k=1$ and renamed 
parameters. From the RTT relations we obtain the following: 
\begin{eqnarray}
%  \label{eq:RTTh23}
  &&     a c x_1 + c a x_2 + c^2 x_3 = 0 \nn & &
       a c - c a + c^2 x_2 = 0 \nn & & 
        - a c + c a + c^2 x_1 = 0 \nn & & 
        - c^2 x_1 + c d - d c = 0 \nn & & 
        - c^2 x_2 - c d + d c = 0 \nn & & 
        c^2 x_3 + c d x_1 + d c x_2 = 0 \nn && 
       a d - c a x_1 + c d x_2 - d a = 0 \nn &&
 - a c x_2 - a d + d a + d c x_1 = 0 \nn & & 
       b c - c a x_2 - c b + d c x_2 = 0 \nn & & 
        - a c x_1 - b c + c b + c d x_1 = 0 \nonumber 
\end{eqnarray} 
\begin{eqnarray}
  \label{eq:RTTh23}
 && 
        - a^2 x_1 + a b + a d x_1 - b a + c b x_2 + c d x_3 = 0 \nn & & 
        - a^2 x_2 - a b + b a + b c x_1 + d a x_2 + d c x_3 = 0 \nn & & 
       b d - c a x_3 - c b x_1 - d a x_2 - d b + d^2 x_2 = 0 \nn & & 
        - a c x_3 - a d x_1 - b c x_2 - b d + d b + d^2 x_1 = 0 \nn & & 
        - a^2 x_3 - a b x_1 - b a x_2 + b d x_1 + d b x_2 + d^2 x_3 = 0 
\end{eqnarray} 
It is necessary to consider several cases.\\
1. In the case $x_1+x_2\neq 0$ (and arbitrary $x_3$) 
from the above follow:
\begin{eqnarray}
  \label{eq:RTTh23a}
  &&c^2 = 0, \qquad   ca \ =\  ac\ =\ 0, 
  \qquad    dc \ =\  cd\ =\ 0, \nonumber\\
  && da \ =\  ad, \qquad cb \ =\ bc, \nonumber\\
  &&   a^2 \ =\ d^2\ =\ ad+bc \nn  
  && ab \ =\   bd \ = \ ba +   (x_1-x_2) bc, 
  \qquad
  db \ =\  bd +    (x_2-x_1) bc 
\end{eqnarray} 
These relations make the resulting algebra rather degenerate. 
Moreover, in order to build a PBW basis we have to look also for 
higher order relations. For instance, using these relations 
we obtain:
$$ a^3 = a^2 d + a b c = a^2 d $$
Furthermore one can eliminate $bc= a^2-ad$, and 
$bd = ba + (x_1-x_2) (a^2-ad)$. 
{}From all these follows that the PBW basis may have only the 
following monomials:
\begin{equation}
  \label{eq:PBWa}
b^n a^k  \,, \quad a^\ell d\,, \quad c\,, \qquad 
n,k\in \ZZ_+ \,, \quad \ell=0,1.
\end{equation}
2. Next we consider the case $x_1 = -x_2 = h$, $x_3\neq -h^2$
then from (\ref{eq:RTTh23}) follow:
\begin{eqnarray}
  \label{eq:RTTh23b}
&&c^2 = 0, \qquad   ca \ =\  ac\ =\ 0, 
\qquad    dc \ =\  cd\ =\ 0, \nonumber\\
  && da = ad, \qquad cb \ =\ bc, \qquad   a^2 \ =\ d^2 \nn  
 && ab = ba + h(a^2 +bc-ad), \qquad
  db = bd - h(a^2 +bc-ad). 
\end{eqnarray}
 The resulting algebra is also degenerate, though 
the above relations are less restrictive than the 
previous case and the possible PBW basis is richer:
\begin{equation}
  \label{eq:PBWb}
b^n a^k d^\ell  \,, \quad b^n c^\ell\,, \qquad 
n,k\in \ZZ_+ \,, \quad \ell=0,1.
\end{equation}
3. Finally in the case  $x_1 = -x_2 = h$, $x_3= -h^2$ this 
coincides with the case $R_{H1,3}$ when $g=-h$. 

 The fifth case $R_{S0,1}$ is a special case of $R_{H2,3}$ when 
$x_1=x_2=0$, $x_3=1$. Thus, the resulting algebra relations 
are obtained from (\ref{eq:RTTh23b}) setting $h=0$.

In the sixth case $R_{S0,2}$ the RTT relations give:
%% \begin{eqnarray}
%%   \label{eq:RTTs02}
%%       c^2 = 0 \nn
%%       a c +c a = 0 \nn 
%%       c d + d c = 0 \nn
%%        - a d + d a = 0\nn
%%        - b c + c b = 0 \nn
%%       a b + b a + c d = 0 \nn
%%       a b + b a + d c = 0 \nn
%% b d + d b + c a = 0 \nn
%% b d + d b + a c = 0 \nn     
%%   a^2 - d^2 = 0 
%% \end{eqnarray}
%% {}From these we obtain the following relations:
\begin{eqnarray}
  \label{eq:RTTs02a}
&&c^2\ =\ 0, \qquad   ca \ =\  ac\ =\ 0, 
\qquad    dc \ =\  cd\ =\ 0, \nonumber\\
  && da\ =\ ad, \qquad cb \ =\ bc, \qquad   a^2 \ =\ d^2  \nn  
 && ab +ba\ =\ 0 , \qquad
  db + bd\ =\ 0 
\end{eqnarray} 
This is a superalgebra, also degenerate like the 
previous two cases. The PBW basis would be as in 
(\ref{eq:PBWb}).

Finally, we are left with the seventh case, which we have anticipated 
to be a hybrid one by the notation $R_{q,h}$. Note that 
setting $q=1$, 
\begin{equation}
  \label{eq:R1h}
  R_h = R_{1,h} = \left(
    \begin{array}{cccc}
      1 & 0 & 0 & h \\
      0 & 1 & 0 & 0 \\
      0 & 0 & 1 & 0 \\
      0 & 0 & 0 & -1   
    \end{array}
  \right) 
\end{equation}
still depends on $h$ and is triangular (in the sense
$R_{21}=R^{-1}$). 

Now we obtain the RTT relations by implementing
(\ref{eq:Rqh}) in (\ref{eq:RTT}) with (\ref{eq:T}):
\begin{eqnarray}
  \label{eq:RTTgl11}
  ba &=& ab+hcd, \qquad\qquad\qquad
  ca \ =\  q^{-1}ac, \nonumber\\
  da &=& ad-(1-q^{-1})bc, \qquad
  cb \ =\  q^{-1}bc, \nonumber\\
  db &=& -bd+hq^{-1}ac, \qquad\qquad
  dc \ =\  -qcd, \nonumber\\
  c^2 &=& 0, \qquad\qquad\qquad\qquad
  ha^2 \ =\  (q+1)b^2+hd^2.
\end{eqnarray} 
This is a superalgebra which we shall denote by $GL_{q,h}(1|1)$. It 
was first written in \cite{FHRclass}, 
where also the dual quantum algebra was given. This is indeed a 
hybrid deformation of $GL(1|1)$ since 
it is known from \cite{Hietarinta} that no transformation of the form
\begin{equation}
  \label{eq:transf}
  R_{qh}\longrightarrow(M\otimes M)R_{qh} (M^{-1}\otimes M^{-1})
\end{equation}
can lead to $R_q$. 

Using the transformation (\ref{eq:transf}) with 
\begin{equation}
  M = \left(
    \begin{array}{cc}
      x & y \\
      0 & x^{-1}
    \end{array}
  \right)
\end{equation}
one gets
\begin{equation}
  \label{eq:Rqhtransf}
  R_{q,h} = \left(
    \begin{array}{cccc}
      1 & 0 & 0 & hx^4-(1+q)(xy)^2 \\
      0 & 1 & 1-q & -2xy \\
      0 & 0 & q & -2qxy \\
      0 & 0 & 0 & -q   
    \end{array}
  \right) 
  \;.
\end{equation}
For no choice of $x,y$ one obtains for the transformed $R$ the form
$R_{q;h=0}$, in contrast with the results of Sect. \ref{sect:Kup} and
\ref{sect:BHP}. 

For $1+q\neq 0$ one may choose to eliminate the top right hand element
by setting $y=\pm x \left(\frac{h}{1+q}\right)^{1/2}$. Denoting
$h'=-2xy$, 
\begin{equation}
  \label{eq:Rqhtransf2}
  R_{q,h} = \left(
    \begin{array}{cccc}
      1 & 0 & 0 & 0 \\
      0 & 1 & 1-q & h' \\
      0 & 0 & q & qh' \\
      0 & 0 & 0 & -q   
    \end{array}
  \right) 
  \;.
\end{equation}

In conclusion we repeat that in order to display the sharp contrast
between the mixed deformations of Sect. \ref{sect:Kup} and
\ref{sect:BHP} and the present case we have restricted our
considerations in both cases to similarity transformations,  
i.e. coboundary 
twists $\cF\equiv (g^{-1}\otimes g^{-1}) \Delta(g)$, which
transform $\cR$ matrices as 
$\cR^\cF \equiv \cF_{21} \cR \cF^{-1} = (g^{-1}\otimes g^{-1}) \cR
(g\otimes g)$.  
Other interesting aspects can be explored by implementing
more general twists. We refer to the discussion in \cite{KulishSym}
concerning the three deformed versions of $gl(1|1)$.

\paragraph{Acknowledgments:} One of us (A.C.) acknowledges with
  pleasure interesting discussions with Petr Kulish. 
~V.K.D. would like to thank for hospitality CERN-TH, where part 
of the work was done. 
~This work was supported in part by the CNRS-BAS France/Bulgaria
  agreement    number 6608.

\end{document}